\documentclass[lettersize,journal]{IEEEtran}
\usepackage{cite}
\usepackage{amsmath,amssymb,amsfonts}
\usepackage{graphicx}
\usepackage{textcomp}
\usepackage{float}
\usepackage{multicol}
\usepackage{amsmath}
\usepackage{amsthm}
\usepackage{amssymb}
\usepackage{amsfonts}
\usepackage{graphicx}
\usepackage{subfigure}
\usepackage{url}
\usepackage{booktabs}
\usepackage{float}
\usepackage{color}
\usepackage{bm}
\usepackage{algorithm}
\usepackage{algpseudocode}

\newtheorem{assmp}{\bf Assumption}

\newtheorem{rem}{\bf Remark}

\newtheorem{lem}{\bf Lemma}

\newtheorem{thm}{\bf Theorem}

\newtheorem{problem}{\bf Problem}

\newtheorem{prop}{\bf Proposition}

\newcommand{\EQ}{\begin{eqnarray}}
	\newcommand{\EN}{\end{eqnarray}}
\newcommand{\EQQ}{\begin{eqnarray*}}
	\newcommand{\ENN}{\end{eqnarray*}}

\ifCLASSINFOpdf
\else
\fi

\begin{document}
	\title{A Refined Algorithm for the Adaptive Optimal Output Regulation Problem}
	
	\author{Liquan~Lin~and~Jie~Huang,~\IEEEmembership{Fellow,~IEEE}
		\thanks{This work was supported in part by the Research Grants Council of the Hong Kong Special Administration Region under grant No. 14201420 and in part by National Natural Science Foundation of China under Project 61973260.}
\thanks{The authors are with the Department of Mechanical and Automation Engineering, The Chinese University of Hong Kong, Hong Kong (e-mail: lqlin@mae.cuhk.edu.hk; jhuang@mae.cuhk.edu.hk. Corresponding author: Jie Huang.)}}
	

	\maketitle
	
	\begin{abstract}
		Given a linear unknown system with  $m$ inputs, $p$ outputs, $n$ dimensional state vector, and $q$ dimensional exosystem,
		the problem of the adaptive optimal output regulation of this system boils down to iteratively solving a set of linear equations and each of these equations contains $\frac{n (n+1)}{2} + (m+q)n$ unknown variables. In this paper, we {refine} the existing algorithm by
		decoupling each of these linear equations into two lower-dimensional linear equations. The first one contains $nq$ unknown variables,  and the second one contains  $\frac{n (n+1)}{2} + mn$ unknown variables.
		As a result, the solvability conditions for these equations are also significantly weakened.
	\end{abstract}
	

	\IEEEpeerreviewmaketitle	
	
	\section{Introduction}
	\IEEEPARstart{R}ecently, reference \cite{gao2016adaptive} studied the optimal output regulation problem for linear unknown systems by state feedback control. The approach in \cite{gao2016adaptive} boils down to iteratively solving a set of linear equations.  For a linear system with  $m$ inputs, $p$ outputs, $n$ dimensional state vector, and $q$ dimensional exosystem, each of these equations contains $\frac{n (n+1)}{2} + (m+q)n$ unknown variables.

	The main objective of  this paper is to refine the  algorithm of \cite{gao2016adaptive}.
	We first show that the $\frac{n (n+1)}{2} + (m+q)n$ unknown variables governed by each linear equation in \cite{gao2016adaptive} can be separated into two groups.
	The first group consists of $nq$ unknown  variables and the second group consists of $\frac{n (n+1)}{2}+mn$ unknown variables.
	As a result, the refined algorithm significantly reduces the computational complexity of the algorithm in \cite{gao2016adaptive}.
	Moreover, the success of the algorithm in \cite{gao2016adaptive} critically hinges on the satisfaction of $(n-p)q+2$ rank conditions for linear equations with  $\frac{n (n+1)}{2} + (m+q)n$ columns.
	In contrast, the refined algorithm only needs to satisfy one rank condition for one linear equation with $\frac{n (n+1)}{2} + (m+q)n$ columns and $(n-p)q+1$ rank conditions of lower dimensional matrices.

It is noted that the refined algorithm also applies to the  adaptive cooperative optimal output regulation of linear multi-agent unknown systems developed in \cite{gao2017cooperative }.

	\indent \textbf{Notation} Throughout this paper, $\mathbb{R}, \mathbb{N} , \mathbb{N}_+$ and $ \mathbb{C}_-$ represent the sets of real numbers, nonnegative integers,  positive integers and the open left-half complex plane, respectively.
$|\cdot|$ represents the Euclidean norm for vectors and the induced norm for matrices.  For $b=[b_1, b_2, \cdots, b_n]^T\in \mathbb{R}^n$, $\text{vecv}(b)=[b_1^2,b_1b_2,\cdots,b_1b_n,b_2^2,b_2b_3,\cdots, b_{n-1}b_n,b_n^2]^T \in \mathbb{R}^{\frac{n(n+1)}{2}}$. For a symmetric matrix $P=[p_{ij}]_{n\times n}\in \mathbb{R}^{n\times n}$, $\text{vecs}(P)=[p_{11},2p_{12},\cdots,2p_{1n},p_{22},2p_{23},\cdots, 2p_{n-1,n}, p_{nn}]^T\in \mathbb{R}^{\frac{n(n+1)}{2}}$. For  $v\in \mathbb{R}^n$, $|v|_P=v^TPv$. For column vectors $a_i, i=1,\cdots,s$,  $\mbox{col} (a_1,\cdots,a_s )= [a_1^T,\cdots,a_s^T  ]^T,$ and, if  $A = (a_1,\cdots,a_s )$, then  vec$(A)=\mbox{col} (a_1,\cdots,a_s )$.
For $A\in \mathbb{R}^{n\times n}$, $\sigma(A) $ denotes the spectrum of $A$, and Tr$(A)$ the trace of $A$. 'bolckdiag' denotes the block diagonal matrix operator. 
	
	\section{Preliminaries}\label{sec2}
	\indent In this section, we review some existing results on the adaptive optimal output regulation problem based on  \cite{gao2016adaptive}  \cite{su2011cooperative}    \cite{huang2004nonlinear} 	 \cite{Kleinman}.
	\subsection{Output Regulation Problem (ORP)}
	Consider a continuous-time linear system in the following form:
	\begin{align} \label{single}
\begin{split}
\dot{x}&=Ax+Bu+Dv\\
			e&=Cx+Fv
		\end{split}
\end{align}
where $x\in \mathbb{R}^{n}, u\in \mathbb{R}^{m}$ and $e\in \mathbb{R}^{p}$ are the state, control input and tracking error of the system, and $v\in \mathbb{R}^{q}$ is the state generated by the following exosystem:
\begin{align}
	\dot{v}&=Ev \label{leaders}	
	\end{align}
Without loss of generality, we assume $C$ is of full row rank.


\begin{problem}\label{outregu}
{\bf [The State Feedback Output Regulation Problem]}
Design a state feedback control law of the following form:
	\begin{align}\label{controllersin}
			u=-Kx+Lv
		\end{align}
where $K$ is called the feedback gain matrix and $L$ is called the feedforward gain matrix  such that the closed-loop system is exponentially stable in the sense that $\sigma(A-BK)\subset \mathbb{C}_-$ and the tracking error $e$ satisfies $\lim\limits_{t\to \infty}e=0$.
\end{problem}

The solvability of Problem \ref{outregu} entails the following assumptions.
    \begin{assmp} \label{ass5}
    	$(A,B)$ is stabilizable.
    \end{assmp}
    \begin{assmp}\label{ass6} There exist a pair of the  matrices
    	$(X,U)$ that satisfies  the following so-called regulator equations:
    	\begin{align}\label{regusin1}
    		XE=&AX+BU+D\\ \label{regusin2}
    		0=&CX+F
    	\end{align}
   \end{assmp}

    The following result gives the solvability of Problem \ref{outregu} \cite{huang2004nonlinear}:
    \begin{thm}
    	Under Assumptions \ref{ass5} and \ref{ass6}, there exist feedback gains $K$ and the feedforward gains  $L=U+KX$ such that (\ref{controllersin}) solves Problem \ref{outregu}.
    \end{thm}

To see the role of Assumption \ref{ass6}, let $\bar{x}:=x -X^*v, \bar{u}:=u-U^*v$, where the pair $(X^*,U^*)$ is any one of the  solutions to  \eqref{regusin1} and
\eqref{regusin2}. Then, we obtain the error system as follows:
    \begin{align}\label{errsyssin}
    	\dot{\bar{x}}&=A\bar{x}+B\bar{u}\\
    	e&=C\bar{x}
    \end{align}
Thus, as long as the state feedback control $\bar{u} = - K \bar{x}$ stabilizes \eqref{errsyssin}, then $\lim_{t \rightarrow \infty} e (t) = 0$.
That is, $u = \bar{u} + U^*v = - K (x  -X^*v) + U^*v = - K x + L v$ where $L = U^*+ K X^*$ solves the ORP.

Under Assumption \ref{ass5},  the following algebraic Riccati equation
    	\begin{align}\label{ARE}
    		A^TP^*+P^*A+C^TQC-P^*BR^{-1}B^TP^*=0
    	\end{align}
has a unique positive definite solution $P^*$. Then, a particular feedback gain $K^*$ is given by $K^*=R^{-1}B^TP^*$.

This $K^*$ is such that the  state feedback controller $\bar{u}^*=-K^*\bar{x}$
solves the following  LQR  problem:
\begin{align}
    		\begin{split}\label{costp2sin}
    			&\min_{\bar{u}}\int_{0}^{\infty}(|e|_{Q}+|\bar{u}|_{R})dt\\
    			&\textup{subject to } (\ref{errsyssin})
    		\end{split}		
    	\end{align}
    	where $Q=Q^T\geq0, R=R^T>0$, with $(A, \sqrt{Q}C)$ observable.

An  algorithm to obtain $P^* $ is given
by iteratively solving the following linear Lyapunov equations  \cite{Kleinman}:
   \begin{align} \label{aleq1sin}
    			\begin{split}
    				&0=(A-BK_{j})^TP_{j} +P_{j}(A-BK_{j})+C^TQC+K_{j}^TRK_{j}
    			\end{split}\\ \label{aleq2sin}
    			&K_{j+1}=R^{-1} B^TP_{j}
    \end{align}
where $j=0,1,\cdots$, and $K_{0}$ is such that $A-BK_{0} $ is a Hurwitz matrix. The algorithm \eqref{aleq1sin} and \eqref{aleq2sin}  guarantees the following properties for $j\in \mathbb{N}$:
    	\begin{enumerate}
    		\item $\sigma(A-BK_{j}) \subset \mathbb{C}_-$;
    		\item $P^*\leq P_{j+1}\leq P_j$;
    		\item $\lim\limits_{j\to \infty}K_j=K^*,\lim\limits_{j\to \infty}P_j=P^*$.
   \end{enumerate}

  \subsection{Adaptive Optimal Output Regulation Problem (AOORP)}

When the matrices $A,B$ and $D$ are unknown, reference \cite{gao2016adaptive} proposed  the integral reinforcement learning (IRL) method to solve the ORP and the LQR problem \eqref{costp2sin}.
This subsection summarizes  the approach in \cite{gao2016adaptive} as follows.

	Given constant matrices $X_{0},X_{1}$ such that $X_{0}=0_{n\times q}, CX_{1}+F=0$, let $h=(n-p)q$ and
 select matrices $X_{i}, i=2,3,\cdots, h+1$, such that all the vectors $\text{vec}(X_{i})$ form a basis for $\text{ker}(I_{q}\otimes C)$ for the positive integer $h$. Then a representation of the general solution to (\ref{regusin2}) is  as follows:
 \begin{align}\label{eq1}
		X&=X_{1}+\sum_{i=2}^{h+1}\alpha_{i}X_{i}
	\end{align}
where $\alpha_{2}, \cdots,  \alpha_{h+1}\in \mathbb{R}$.

Define a Sylvester map $\mathcal{S}(X_i)=X_{i}E-AX_{i}$. Then,
(\ref{regusin1})  can be put as follows:
	\begin{align}\label{eq2}
		\mathcal{S}(X)&=\mathcal{S}(X_{1})+\sum_{i=2}^{h+1}\alpha_{i}\mathcal{S}(X_{i})=BU+D
			\end{align}

	Combining (\ref{aleq2sin}), (\ref{eq1}) and (\ref{eq2}) gives
	\begin{align} \label{regunew}
		\mathcal{A}\chi=b
	\end{align}
	where
	\begin{align*}
		\mathcal{A}&=\begin{bmatrix}
			\mathcal{A}_{1} & \mathcal{A}_{2}
		\end{bmatrix}\\
		\mathcal{A}_{1}&=\begin{bmatrix}
			\text{vec}(\mathcal{S}(X_{2})) & \cdots & \text{vec}(\mathcal{S}(X_{h+1}))\\
			\text{vec}(X_{2}) & \cdots & \text{vec}(X_{h+1})
		\end{bmatrix}\\
		\mathcal{A}_{2}&= \begin{bmatrix}
			0 & -I_q \otimes (P_{j}^{-1}K_{j+1}^TR)\\
			-I_{nq} & 0
		\end{bmatrix}\\
		\chi&=[\alpha_{2}, \cdots , \alpha_{h+1}, \text{vec}(X)^T, \text{vec}(U)^T]^T\\
		b&=\begin{bmatrix}
			\text{vec}(D-\mathcal{S}(X_{1}))\\
			-\text{vec}(X_{1})
		\end{bmatrix}	
	\end{align*}
	\begin{rem}
		Equation (\ref{regunew}) indicates that, with the  information of $X_i,\mathcal{S}(X_{i}),\; i=1,\cdots,h+1$, $D$, $P_j$ and $K_{j+1}$ for some $j$, one can
obtain the solution  to the regulator equations (\ref{regusin1}) and (\ref{regusin2}) without knowing the matrices $A$ and $B$.
	\end{rem}

Let $A_{j}=A-BK_{j}$ with $K_{j}$ being obtained from \eqref{aleq2sin}, and for $i=1,\cdots,h+1$,
 let $\bar{x}_{i}=x-X_{i}v$ for $i=0,1,\cdots, h+1$. Then the dynamics of the error systems can be derived as follows:
		\begin{align}\label{errsin}
		\begin{split}
			\dot{\bar{x}}_{i}=&Ax+Bu+Dv-X_{i}Ev\\
			=&A_{j}\bar{x}_{i}+B(K_{j}\bar{x}_{i}+u)+(D-\mathcal{S}(X_{i}))v
		\end{split}
	\end{align}
	
	
For any $ t \geq 0$, $\delta t >0$, using (\ref{aleq1sin}), (\ref{aleq2sin}) and \eqref{errsin} gives	
	\begin{align}
		\begin{split} \label{IRLsin1}
			&|\bar{x}_{i}(t+\delta t)|_{P_{j}}-|\bar{x}_{i}(t)|_{P_{j}}\\
			=& \int_{t}^{t+\delta t} [|\bar{x}_{i}|_{A_{j}^TP_{j}+P_{j}A_{j}}+2(u+K_{j}\bar{x}_{i})^TB^TP_{j}\bar{x}_{i}\\ & +2v^T(D-\mathcal{S}(X_{i}))^TP_{j}\bar{x}_{i}]d\tau\\
			=& \int_{t}^{t+\delta t} [-|\bar{x}_{i}|_{C^TQC+K_{j}^TRK_{j}}+2(u+K_{j}\bar{x}_{i})^TRK_{j+1}\bar{x}_{i}\\ &+2v^T(D-\mathcal{S}(X_{i}))^TP_{j}\bar{x}_{i}]d\tau\\
		\end{split}
	\end{align}

For any vectors $a\in \mathbb{R}^n,b \in \mathbb{R}^m$ and any integer $s\in \mathbb{N}_+$, define
	
	\begin{align}
		\begin{split}\label{newdefi}
			\delta _a=&[\text{vecv}(a(t_1))-\text{vecv}(a(t_0)), \cdots ,\\ &\text{vecv}(a(t_s))-\text{vecv}(a(t_{s-1}))]^T\\
			\Gamma_{ab}=&[\int_{t_0}^{t_1}a\otimes b d\tau , \int_{t_1}^{t_2}a\otimes b d\tau, \cdots , \int_{t_{s-1}}^{t_s}a\otimes b d\tau]^T
		\end{split}
	\end{align}
	
	Then, using the notation (\ref{newdefi}), equation (\ref{IRLsin1}) can be arranged as the following linear equation:
	
	\begin{align}\label{IRLsinlinear}
		\Psi_{ij}
		\begin{bmatrix}
			\text{vecs} (P_{j})\\
			\text{vec}(K_{j+1})\\
			\text{vec}((D-\mathcal{S}(X_{i}))^TP_{j})
		\end{bmatrix}=\Phi_{ij}
	\end{align}
	where
	\begin{align}
		\Psi_{ij}&=[\delta_{\bar{x}_{i}}, -2\Gamma_{\bar{x}_{i}\bar{x}_{i}}(I_{n}\otimes K_{j}^TR)-2\Gamma_{\bar{x}_{i}u}(I_{n}\otimes R), -2\Gamma_{\bar{x}_{i}v}] \notag \\
		\Phi_{ij}&=-\Gamma_{\bar{x}_{i}\bar{x}_{i}}\text{vec} (C^TQC+K_{j}^TRK_{j})\label{psijk}
	\end{align}
	
	The following lemma \cite{gao2016adaptive} gives the condition which guarantees the uniqueness of the solution to (\ref{IRLsinlinear}).
	\begin{lem}\label{ranconsinlem}
		For  $i\in \{0,1,2,\cdots, h+1\}$, if there exist a $s^*\in \mathbb{N}_+$ such that for all $s>s^*$, and for all sequences $t_0<t_1<\cdots <t_s$,
		\begin{align} \label{rankconsinpre}
			\textup{rank} ([\Gamma_{\bar{x}_{i}\bar{x}_{i}}, \Gamma_{\bar{x}_{i}u}, \Gamma_{\bar{x}_{i}v}])= \frac{n(n+1)}{2}+(m+q)n
		\end{align}
		then the matrix $\Psi_{ij}$ has full column rank for all $j\in \mathbb{N}$.
	\end{lem}

	\section{Refinements of the Algorithm for AOORP} \label{sec3}
  In this section, we refine the algorithm for solving  AOORP.

Let us first give a specific way to find the null space of $C X=0$ where
$X \in \mathbb{R}^{n \times q}$ as follows: 	
	
\begin{prop}
	Let $\{y_1,y_2,\cdots, y_{n-p}\}$ form a basis for the null space of $Cy=0$ with $y \in \mathbb{R}^{n}$. Then the following matrices form a basis for the null space of $C X=0$ where
$X \in \mathbb{R}^{n \times q}$:
\begin{align}\label{Xi}
			X_{(n-p)k+i}=[0_{n\times k}, y_i, 0_{n\times (q-k-1)}]
		\end{align}
	for $k=0,\cdots, q-1 \text{ and } i=1,\cdots,n-p$.
\end{prop} 	

Next we show that (\ref{IRLsinlinear}) can be reduced to two  lower-dimensional linear equations.
	
	\begin{enumerate}
		\item Obtain $D, P_{0}, K_{1}$ by solving the following equation
\begin{align}\label{IRLsinlinear0}
		\Psi_{00}
		\begin{bmatrix}
			\text{vecs} (P_{0})\\
			\text{vec}(K_{1})\\
			\text{vec}(D^T P_{0})
		\end{bmatrix}=\Phi_{00}
	\end{align}
which is obtained from
 (\ref{IRLsinlinear}) by letting $j=0,i=0$.
By Lemma \ref{ranconsinlem} with $i=0$, \eqref{IRLsinlinear0} is solvable if
\begin{align} \label{rankconmul10}
		\textup{rank} ([\Gamma_{\bar{x}_{0}\bar{x}_{0}}, \Gamma_{\bar{x}_{0}u}, \Gamma_{\bar{x}_{0}v}])= \frac{n(n+1)}{2}+(m+q)n
	\end{align}
		\item Let $\Psi_{ij}=[\Psi_{ij}^1, \Psi_{ij}^2, \Psi_{ij}^3]$ with $\Psi_{ij}^1=\delta_{\bar{x}_{i}}, \Psi_{ij}^2=-2\Gamma_{\bar{x}_{i}\bar{x}_{i}}(I_{n}\otimes K_{j}^TR)-2\Gamma_{\bar{x}_{i}u}(I_{n}\otimes R), \Psi_{ij}^3=-2\Gamma_{\bar{x}_{i}v}$. Then,
		from  \eqref{IRLsinlinear} with  $j=0$, we have
		\begin{align*}
			& \Psi_{i0}^3 \text{vec}((D-\mathcal{S}(X_{i}))^TP_{0}) \notag \\
			=& -\Psi_{i0}^1\text{vecs} (P_{0})- \Psi_{i0}^2\text{vec}(K_{1}) + \Phi_{i0}
		\end{align*}
which can be put in the following form:
\begin{align}  \label{new2}
			& \Gamma_{\bar{x}_{i}v}  \text{vec}((D-\mathcal{S}(X_{i}))^TP_{0}) \notag \\
			=& \frac{1}{2} (\Psi_{i0}^1\text{vecs} (P_{0})+\Psi_{i0}^2\text{vec}(K_{1})- \Phi_{i0})
		\end{align}
If
\begin{align}\label{rankconsin2}
		\textup{rank}(\Gamma_{\bar{x}_{i}v})=qn, \forall i=1,2,\cdots, h+1,
	\end{align}
then we can solve $\mathcal{S}(X_{i})$ for $i = 1, \cdots, h+1$ from  \eqref{new2}.
		\item Let $M\in R^{(n^2)\times(\frac{n\times (n+1)}{2})}$ be a constant matrix such that $M\text{vecs}(P_{j})=\text{vec}(P_{j})$ and $M_{D}=(I_{n}\otimes D^T) M$. Then, we have
		\begin{align}\label{new1}
			&\text{vec}(D^TP_{j}) \notag \\=&(I_{n}\otimes D^T )\text{vec}(P_{j}) \notag \\
			=&(I_{n}\otimes D^T )M\text{vecs}(P_{j}) \notag \\
			=&M_{D}\text{vecs}(P_{j})
		\end{align}
		As a result,    \eqref{IRLsinlinear} with $i=0$ reduces to the following:
		\begin{align} \label{new3}
			&\Psi_{0j}
			\begin{bmatrix}
				\text{vecs} (P_{j}) \\
				\text{vec}(K_{j+1})\\
				\text{vec} (D^T P_{j})
			\end{bmatrix} \notag  \\
			=&\Psi_{0j}^1\text{vecs} (P_{j})+\Psi_{0j}^2\text{vec}(K_{j+1}) \notag  \\&+\Psi_{0j}^3\text{vec} (D^T P_{j}) \notag  \\
			=&\Psi_{0j}^1\text{vecs} (P_{j})+\Psi_{0j}^2\text{vec}(K_{j+1}) \notag  \\&+\Psi_{0j}^3M_{D}\text{vecs}(P_{j}) \notag  \\
			=& [\Psi_{0j}^1+\Psi_{0j}^3M_{D}, \Psi_{0j}^2]\begin{bmatrix}
				\text{vecs} (P_{j})\\
				\text{vec}(K_{j+1})
			\end{bmatrix} =\Phi_{0j}
		\end{align}
	\end{enumerate}
	
The solvability of (\ref{new3}) is given as follows.

\begin{lem}
		If there exist a $s^*\in \mathbb{N}_+$ such that for all $s>s^*$, and for all sequences $t_0<t_1<\cdots <t_s$,
	\begin{align} \label{rankconsin}
		\textup{rank} ([\Gamma_{\bar{x}_{0}\bar{x}_{0}}, \Gamma_{\bar{x}_{0}u}])= \frac{n(n+1)}{2}+mn
	\end{align}	
		then the matrix $[\Psi_{0j}^1+\Psi_{0j}^3M_{D}, \Psi_{0j}^2]$ has full column rank for all $j\in \mathbb{N}$.
	\end{lem}

	\begin{rem} \label{rem3}
The original algorithm needs to solve (\ref{IRLsinlinear}) for $D, P_j$ and $K_{j+1}$ with $i = 0$ and $j=0, \cdots, j^*$ where $j^*$ is such that $||P_{j^*} - P_{j^*-1}||$ is sufficiently small, and then solve (\ref{IRLsinlinear}) for ${S}(X_i)$ with $j =j^*$ and  $i=1,2,\cdots,h+1$. In particular, it needs the rank condition \eqref{rankconsinpre} to be satisfied for
$i\in \{0,1,2,\cdots, h+1\}$. In contrast, in the refined algorithm,
we only need  the rank condition \eqref{rankconmul10} which is obtained from
the rank condition \eqref{rankconsinpre} with $i=0$ and
the rank conditions (\ref{rankconsin2}) with $i=1,2,\cdots,h+1$ to be satisfied.
 \end{rem}

%

\end{document}